\documentclass[11pt,a4paper]{article}

\usepackage[latin1]{inputenc}
\usepackage{latexsym,amssymb}
\usepackage{amsmath,amssymb,amsfonts,amscd,theorem}
\usepackage{enumerate}
\usepackage[francais,english]{babel}
\usepackage{graphics, graphicx, psfrag}

\newcommand{\C}{\mathbb{C}}

\newcommand{\Q}{\mathbb{Q}}
\newcommand{\R}{\mathbb{R}}
\newcommand{\Z}{\mathbb{Z}}
\newenvironment{proof}[1][Preuve]%
{
\begin{trivlist} \item[]  {\em #1.} }%
{\hspace*{\fill} $\Box$
\end{trivlist}}

{\bigl(
\begin{smallmatrix} }%
{
\end{smallmatrix} \bigr)}

\newtheorem{thm}{Th\'eor\`eme}
\newtheorem{prop}[thm]{Proposition}
\newtheorem{lem}[thm]{Lemme}
\newtheorem{dfn}[thm]{D\'efinition}

\newtheorem{rem}[thm]{Remarque}

\title{Cat\'egorification des coefficients de la matrice de la repr\'esentation de Burau}

\author{Abderrahmane BOUCHAIR\\
D\'epartement de Math\'ematiques\\
Universit\'e de Jijel\\
18000 Jijel,
Alg\'erie\\
courriel : abdebouchair@yahoo.com}
\date{}

\begin{document}
\selectlanguage{francais}

\maketitle

\begin{abstract}

Nous cat\'egorifions explicitement les coefficients de la matrice de la repr\'esentation de Burau en utilisant des m\'ethodes g\'eom\'etriques \'el\'ementaires. Nous montrons que cette cat\'egorification est fid\`ele dans le sens o\`u elle d\'etecte la tresse triviale.
\paragraph*{}
Mots cl\'es: Cat\'egorification, groupe de tresses, homologie de Floer, repr\'esentation de Burau.
\medskip
\noindent 

\end{abstract}

\selectlanguage{english}
\begin{abstract}

We categorify the coefficients  of the Burau representation matrix using elementary geometrical methods. We show the faithfulness of this categorification in the sense that it detects the trivial braid.  
 \end{abstract}

\selectlanguage{francais}

\section{Introduction}
\label{sec:intro}
\begin{figure}[htbp] \label{fig:base}
\psfrag{d0}{$ d_{0} $}
\psfrag{p1}{ $p_{1}$}
\psfrag{p2}{$ p_{2} $}
\psfrag{pn}{$ p_{n} $}
\psfrag{alpha0}{$ \alpha_{0} $}
\psfrag{alpha1}{$ \alpha_{1} $}
\psfrag{alphan-1}{$ \alpha_{n-1} $}
\psfrag{beta1}{$\beta_{1} $}
\psfrag{betan-1}{$ \beta_{n-1}$}
\centerline{\includegraphics[width=8cm ,height=6cm]{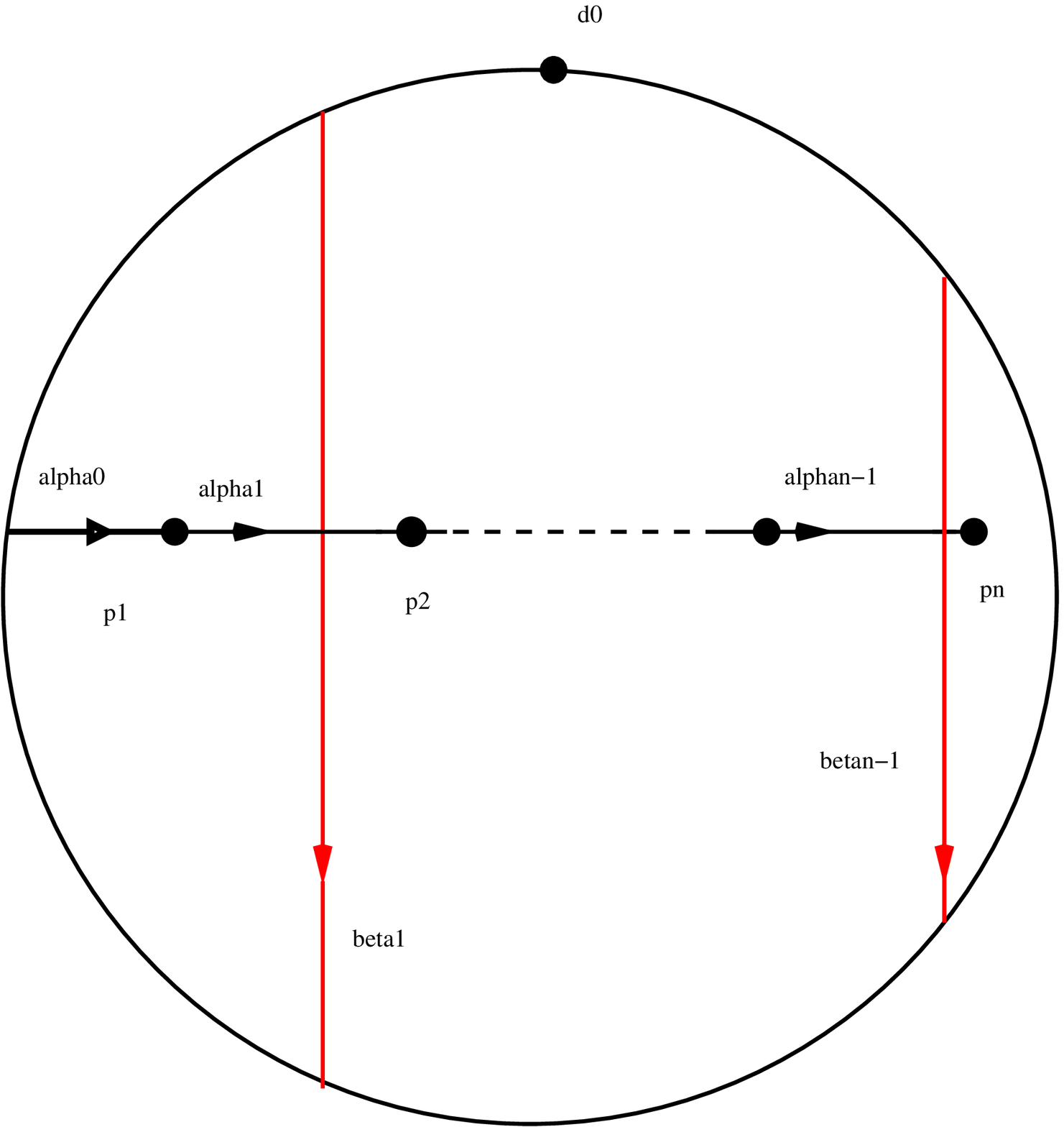}}
\caption{}
\end{figure}

Le groupe de tresse $B_{n}$ introduit par Emile Artin en 1925 est engendr\'e par $ n-1 $ g\'en\'erateurs 
$\sigma_{1}$,...,$\sigma_{n-1}$. Les relations entre les g\'en\'erateurs sont:
\begin{enumerate}
\item $\sigma_i\sigma_j=\sigma_j\sigma_i$   si
$\arrowvert{j-i}\arrowvert\geq2$

\item $\sigma_i\sigma_{i+1}\sigma_i=\sigma_{i+1}\sigma_i\sigma_{i+1}$
 pour  $i=1$,..,$n-2.$
\end{enumerate}
L'une des repr\'esentations classiques des groupes de tresses est la repr\'esentation de Burau, d\'efinie par Burau  en 1936 \cite{BURAU36}, qui envoie le groupe $ B_{n} $ dans le groupe des matrices inversibles de taille $ n\times n $ \`a coefficients dans l'anneau des polyn\^omes de Laurent $ \Z\left[ t,t^{-1}\right]  $ (voir \cite{BIRMAN74} pour plus de d\'etails). 
La repr\'esentation de Burau est r\'eductible, elle se d\'ecompose en une repr\'esentation de dimension $ 1 $ et une repr\'esentation irr\'eductible de dimension $ n-1 $ appel\'ee la repr\'esentation de Burau r\'eduite qu'on  note par $ \widetilde{\rho}:B_n\rightarrow
GL_{n-1}(\Z \left[ t,t^{-1}\right] ), $ d\'efinie comme suit :
\begin{displaymath}
\sigma_i \mapsto I_{i-2}\oplus \left(\begin{array}{ccc} 1&0&0\\
t&-t&1\\ 0&0&1
\end{array}\right) \oplus I_{n-i-2}
\end{displaymath}
O\`u $ -t $ dans le milieu de la matrice de taille $ 3 \times 3 $ est toujours en position $ \left( i,i\right) . $
\paragraph*{}
Cat\'egorifier un invariant qui est un polyn\^ome de Laurent \`a coefficients entiers consiste \`a construire une homologie dont la caract\'eristique d'Euler gradu\'ee est \'egale \`a ce polyn\^ome. L'homologie de Khovanov \cite{KHOVANOV} est une cat\'egorification du polyn\^ome de Jones. Il est naturel de se demander s'il existe une cat\'egorification des coefficients de la matrice $ \widetilde{\rho}(\sigma) ,$ pour une tresse donn\'ee $ \sigma $ dans $ B_{n}. $
\paragraph*{}
En 2001, M. Khovanov et P. Seidel \cite{KH-S} donnent une cat\'egorification de la repr\'esentation de Burau. C'est une action de $ B_{n} $ sur la cat\'egorie deriv\'ee des modules gradu\'es d\'efinis sur une famille de quotients d'alg\`ebres de carquois. Ils expliquent le lien avec le point de vue de la g\'eom\'etrie symplectique. Ils montrent que cette cat\'egorification est fid\`ele.
\paragraph*{}
Dans ce papier, suivant les id\'ees de M. Khovanov et P. Seidel, nous cat\'egorifions, d'une mani\`ere g\'eom\'etrique naive, chaque coefficient de la matrice de la repr\'esentation de Burau et nous montrons le r\'esultat de fid\'elit\'e dans ce cadre. Dans la section $ 2, $ nous donnons une description homologique de la repr\'esentation de Burau. Ensuite, nous d\'efinissons le nombre d'intersection g\'eom\'etrique de deux courbes sur une surface. Soit $ \sigma\in B_{n}, $ soit $ \alpha_{j},\beta_{i} $ deux courbes dans le disque \'epoint\'e $ D_{n} $ comme dans la Figure \ref{fig:base} et $ h_{\sigma} $ un diff\'eomorphisme repr\'esentant la tresse $ \sigma $. Dans la section 3, nous d\'efinissons le complexe de Floer $ \left( CF(\beta_{i},h_{\sigma}(\alpha_{j})),\partial\right)  $ o\`u $CF(\beta_{i},h_{\sigma}(\alpha_{j}))$ est le groupe ab\'elien libre engendr\'e par les points d'intersection de $ \beta_{i} $ avec $ h_{\sigma}(\alpha_{j}). $ Nous montrons que la caract\'eristique d'Euler gradu\'ee de son homologie est \'egale au coefficient $ b_{ij} $ de la matrice $ \widetilde{\rho}(\sigma) $ en position $ (i,j). $ Dans la section 4, nous montrons que cette cohomologie ne d\'epend que de la classe d'isotopie de la tresse i.e si $ h_{\sigma},h'_{\sigma} $ sont deux diff\'eomorphismes (repr\'esentant deux tresses g\'eom\'etriques) isotopes que nous notons par $ h_{\sigma}\simeq h'_{\sigma} $, alors les groupes de cohomologie associ\'es $ HF^{*}(\beta_{i},h_{\sigma}(\alpha_{j})),HF^{*}(\beta_{i},h'_{\sigma}(\alpha_{j}) $ sont isomorphes. Dans la section 5, nous montrons que cette cat\'egorification est fid\`ele dans le sens o\`u elle d\'etecte la tresse triviale.
\section{Pr\'eliminaires} \label{sec:courbe}
\subsection{Repr\'esentation homologique de $ B_{n} $}
Soient $ D =D^{2}$ le disque unit\'e du plan complexe $ \C, $ $ P_{n} $ un ensemble fini de $ n $ points 
$ p_{1}$,...,$p_{n} $ dans l'int\'erieur de $ D $ appel\'es perforations. D\'efinissons $ \text{Diff}^{+}(D ,P_{n}) $ comme \'etant le groupe des diff\'eomorphismes $ h:D\rightarrow D $ qui pr\'eservent l'orientation, laissent le bord fixe point par point et tels que $ h(P_{n})=P_{n}. $ Le groupe de diff\'eotopies "mapping class group" que nous notons par $ \text{MCG}(D , P_{n}) $ est $ \pi_{0}(\text{Diff}^{+}(D ,P_{n})) ,$ i.e le groupe quotient de $ \text{Diff}^{+}(D ,P_{n}) $ par le sous groupe des diff\'eomorphismes de $\text{Diff}^{+}(D ,P_{n}) $ isotopes \`a l'identit\'e relativement \`a  $ \partial D \cup P_{n}.$ Le groupe $ B_{n} $ est isomorphe au groupe $ \text{MCG}(D , P_{n}) .$ 
\paragraph*{}
En utilisant la d\'efinition de $ B_{n} $ comme groupe de diff\'eotopies, la repr\'esentation de Burau peut s'obtenir en faisant agir $ B_{n} $ sur l'homologie du rev\^etement infini cyclique du disque \'epoint\'e. Ce point de vue s'appelle la d\'efinition homologique de la repr\'esentation de Burau (voir \cite{KASSEL-TURAEV} pour plus de d\'etails).
\paragraph*{}
Notons $ D_{n}=D\setminus P_{n}. $ Soit $ d_{0} $ un point base sur le bord de $ D. $ Le groupe $ \pi_{1} \left( D_{n}, d_{0} \right)  $ est libre, engendr\'e par $ n $ g\'en\'erateurs $ x_{1}$,...,$x_{n}, $ o\`u $ x_{i} $ repr\'esente un lacet dans $ D_{n} $ bas\'e en $ d_{0} $ qui tourne dans le sens direct autour du point $ p_{i}. $ Consid\'erons  l'\'epimorphisme de groupe $ \varphi : \pi_{1} \left( D_{n},d_{0} \right) \rightarrow \Z $  d\'efini par
\begin{displaymath}
\gamma = x_{i_{1}}^{n_{1}}x_{i_{2}}^{n_{2}}...x_{i_{r}}^{n_{r}}\longmapsto \varphi \left( \gamma \right) = \sum_{i=1}^{r}n_{i}.
\end{displaymath} 
L'entier $ \varphi\left( \gamma \right)  $ repr\'esente le nombre alg\'ebrique total de tour de $ \gamma $ autour des perforations $ p_{1},p_{2}$,...,$p_{n}. $ Soit $ \tilde{D}_{n} $ le rev\^etement r\'egulier associ\'e au noyau de $ \varphi. $ Le groupe de transformations de rev\^etement de $ \tilde{D}_{n} $ est isomorphe \`a $ \Z=\left\langle t\right\rangle . $ Le groupe d'homologie $ H_{1}( \tilde{D}_{n})$ a une structure de $ \Z \left[ t,t^{-1} \right]-$module libre de rang $ n-1 $, o\`u $ t $ agit par transformation de rev\^etement.
\paragraph*{}
Un diff\'eomorphisme $ h \in \text{Diff}^{+}(D,P_{n})$ ($ h $ repr\'esente un \'el\'ement de $ B_{n} $) se rel\`eve d'une mani\`ere unique \`a un automorphisme $ \tilde{h}: \tilde{D}_{n} \rightarrow \tilde{D}_{n} $ qui fixe la fibre au dessus de $ d_{0} $ point par point. L'automorphisme $\tilde{h}  $ induit un automorphisme de $ \Z \left[ t,t^{-1} \right]-$module $\tilde{h}_{*}:H_{1} ( \tilde{D}_{n} )\rightarrow H_{1} ( \tilde{D}_{n} ).$ L'application $ h \longmapsto \tilde{h}_{\ast} $ d\'efinit un homomorphisme de groupes \'equivalent \`a la repr\'esentation de Burau r\'eduite $ \tilde{\rho}.$
\paragraph*{}
Soit $ \sigma \in B_{n}.$ Supposons que $ \tilde{\rho}(\sigma)=\left( b_{ij}\right)_{1\leq i,j\leq n-1}$ est la matrice de la repr\'esentation de Burau associ\'ee \`a $ \sigma. $ A l'aide de la notion d'intersection alg\'ebrique de deux courbes, nous retrouvons les coefficients $ b_{ij} $ comme suit: Soit $ h_{\sigma}\in \text{Diff}^{+}(D,P_{n}) $ un diff\'eomorphisme (repr\'esentant $ \sigma $). Pour tout $ i=1$,...,$n-1, $ soit $ \alpha_{i} $ le segment de $ p_{i} $ vers $ p_{i+1}, $ et soit $ \beta_{i} $ une corde verticale \`a extr\'emit\'es dans le bord  $ \partial D, $ orient\'e de haut en bas et qui passe entre les perforations $ p_{i} $ et $ p_{i+1}, $ voir la Figure \ref{fig:base}.  Choisissons $ \tilde{d}_{0}\in \partial \tilde{D}_{n} $ un point au dessus de $ d_{0} $ et choisissons deux relev\'es $ \widetilde{h_{\sigma}(\alpha}_{j}), \tilde{\beta}_{i}\subset \tilde{D}_{n} $ de $ h_{\sigma}(\alpha_{j}) $ et $ \beta_{i} $ respectivement. Alors
\begin{eqnarray} \label{eq:Ialgebrique}
b_{ij}=\sum _{k \in \Z}(t^{k}\tilde{\beta}_{i}.\widetilde{h_{\sigma}(\alpha}_{j}))t^{k},
\end{eqnarray}
O\`u $ (t^{k}\tilde{\beta}_{i}.\widetilde{h_{\sigma}(\alpha}_{j})) \in \Z $ est le nombre alg\'ebrique d'intersection standard des arcs orient\'es $ t^{k}\tilde{\beta}_{i} $ et $ \widetilde{h_{\sigma}(\alpha}_{j}) $ dans $ \tilde{D}_{n}. $
\paragraph*{}
Dans tout ce qui suit, $ \sigma \in B_{n}$ est une tresse \`a $ n $ brins et $ h_{\sigma}\in \text{Diff}^{+}(D,P_{n}) $ repr\'esentant $ \sigma. $ Les courbes $ \alpha_{j},\beta_{i} $ sont comme dans la Figure \ref{fig:base} pour certain $ i,j \in \left[ \left[   1,...,n-1 \right] \right]   .$ 
\subsection{Intersection g\'eom\'etrique des courbes}
Soit $ \Sigma $ une surface orientable compacte connexe \'eventuellement \`a bord. Soit $ P_{n} $ un ensemble fini de points dans l'int\'erieur de $ \Sigma. $ Nous utilisons ici la m\^eme notation que Khovanov-Seidel \cite{KH-S}, une courbe de $ \left( \Sigma, P_{n} \right)  $ est un sous ensemble $ c $ dans $ \Sigma $ qui est soit une courbe simple ferm\'ee de $ \Sigma \setminus \left( \partial \Sigma\cup P_{n}\right)  $ et essentielle, soit l'image d'un plongement $\gamma : \left[ 0,1\right]  \rightarrow \Sigma$ transverse au bord telle que $ \gamma^{-1}\left( \partial \Sigma \cup Pn\right) = \left\lbrace  0,1 \right\rbrace . $ Deux courbes $ c_{1},c_{2}$ de $ \left( \Sigma, P_{n} \right)  $ sont isotopes, ce qui se note $ c_{1}\simeq c_{2}, $ s'il existe une isotopie  dans  $\text{Diff}^{+}(\Sigma,P_{n}) $ qui envoie $ c_{1} $ sur $ c_{2}. $ Notons que les extr\'emit\'es sur $ \partial\Sigma $ restent fixe le long de l'isotopie. 
\paragraph*{}
\begin{figure}[htbp] 
\psfrag{c1}{$ c_{1} $}
\psfrag{c2}{$ c_{2} $}
\centerline{\includegraphics[width=3cm,height=2cm]{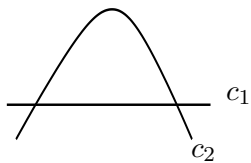}}
\caption{Bigone \'el\'ementaire entre deux courbes}\label{fig:big-elem}
\end{figure}

Soient $ c_{1},c_{2} $ deux courbes simples ferm\'ees et essentielles. Le nombre d'intersection g\'eom\'etrique de $ c_{1} $ et $ c_{2}, $ not\'e $ i(c_{1},c_{2}) ,$ est le nombre minimal de points d'intersection d'un repr\'esentant dans la classe d'isotopie de $ c_{1} $ avec un repr\'esentant de la classe d'isotopie de $ c_{2}. $ Rappelons qu'un bigone (\'el\'ementaire\footnote{Ce mot est ajout\'e pour distinguer entre cette d\'efinition et la d\'efinition \ref{dfn:bigone} qui est plus g\'en\'erale.}) est un disque plong\'e dans $ \Sigma $ dont le bord est form\'e d'un arc de $ c_{1} $ et d'un arc de $ c_{2} $ et son int\'erieur ne coupe pas $ c_{1}\cup c_{2}. $ Nous avons la proposition suivante: 
\begin{prop} \label{prop:bigone}
Soient $ c_{1} $ et $ c_{2} $ deux courbes simples ferm\'ees de $ \Sigma, $ essentielles et se coupant transversalement. Alors $ i \left( c_{1},c_{2} \right) = \vert c_{1}\cap c_{2} \vert $ si et seulement si $ c_{1} $ et $ c_{2} $ ne bordent aucun bigone \'el\'ementaire.
\end{prop}
Pour la preuve voir (\cite{LAUDENBACH}, proposition 3.10) ou (\cite{PARIS}, proposition 3.2).
La preuve de la Proposition \ref{prop:bigone} peut s'adapter au cas des courbes non ferm\'ees $ c_{1},c_{2} $ de $ \left( \Sigma, P_{n} \right)  $ avec $ \partial c_{1}\cap \partial c_{2}=\Phi. $
\paragraph*{}
\begin{figure}[htbp] 
\psfrag{x}{$ x $}
\psfrag{y}{$ y $}
\psfrag{a}{ (a)}
\psfrag{b}{ (b) }
\psfrag{c}{ (c)}
\centerline{\includegraphics[width=8cm, height=4cm]{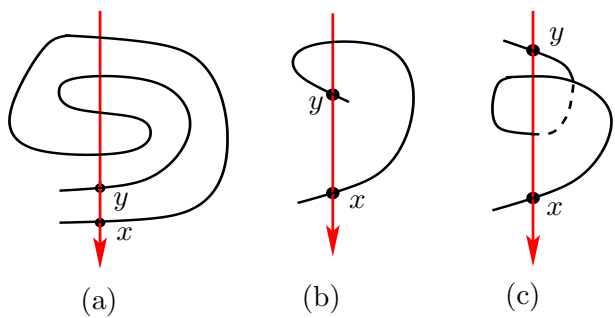}}
\caption{(a) est un bigone de $ x $ \`a $ y $, par contre (b) et (c) ne sont pas des bigones.}\label{fig:immersion}
\end{figure}
Dans \cite{ABOUZAID} M. Abouzaid d\'efinit l'homologie de Floer pour certaine classe de courbes immerg\'ees dans une surface et l'utilise pour d\'efinir la cat\'egorie de Fukaya de cette surface. Nous donnons la d\'efinition suivante:
\begin{dfn} \label{dfn:bigone}
Etant donn\'es deux points $ x , y\in \beta_{i} \cap h_{\sigma}(\alpha_{j}). $ On appelle bigone de $ x $ vers $ y $ une immersion \`a reparam\'etrage orient\'e pr\`es \mbox{$ u: \mathbb{D}\rightarrow D_{n} $} du demi disque \mbox{$ \mathbb{D}:=\{ z\in \C:\Vert z\Vert\leq 1, Re( z)\geq 0\}, $} qui envoie $ -i $ sur $ x, $ $ i $ sur $ y $ et v\'erifie les conditions suivantes:
\begin{enumerate}
\item $ u(\mathbb{D}\cap i\R ) \subset \beta _{i}, $
\item $ u(\mathbb{D}\cap S^{1}) \subset h_{\sigma}(\alpha_{j}). $
\end{enumerate}
\end{dfn}
\begin{figure}[htbp] 
\psfrag{x}{$ x $}
\psfrag{y}{$ y $}
\psfrag{snega}{$ \varepsilon (u)=-1 $}
\psfrag{sposi}{$ \varepsilon (u)=1 $}
\centerline{\includegraphics[width=6cm,height=3cm]{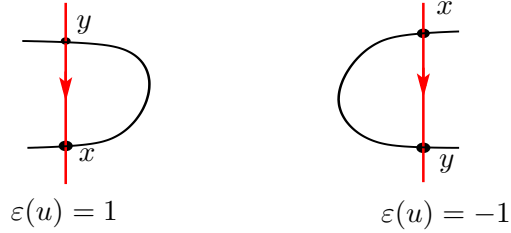}}
\caption{Bigone $ u $ de $ x $ \`a $ y. $}\label{fig:signe-bigone}
\end{figure}
On note $ \mathcal{M} \left( x,y \right)  $ l'ensemble des bigones de $ x $ vers $ y. $ Dans notre cas on prend des immersions, voir la Figure \ref{fig:immersion}. Si $ u\in \mathcal{M} \left( x,y \right),  $
on note $ \alpha_{xy}$ et $ \beta_{xy} $  les sous arcs de $ h_{\sigma}(\alpha_{j}) $ et $ \beta_{i} $ respectivement, qui bordent le bigone $ u. $ 
Dans tout le texte sauf mention du contraire, on suppose que $ \beta_{i} $ est une corde verticale. On associe \`a $ u $ un signe $ \varepsilon ( u) $ d\'efini comme suit: $ \varepsilon (u)=1 $ si l'orientation de $\beta_{i}$ coinc\"ide avec son orientation comme bord de $ u $ et $ \varepsilon  (u)=-1 $ sinon, voir la Figure \ref{fig:signe-bigone}.
\paragraph*{}
Soit $ x $ un point d'intersection, et soit $ V_{x} $ un voisinage de $ x $ dans $ D_{n}. $ Un model local de $ x $ est une application orient\'ee de $ V_{x} $ dans un voisinage de l'origine dans $ \R^{2} $ qui envoie $ h_{\sigma}(\alpha_{j}) $ sur l'axe des abscisses et $ \beta_{i} $ sur l'axe des ordonn\'ees. Cela nous permet de comparer les ordonn\'ees des points d'intersection sur la courbe $ \beta_{i}. $ Rappelons que nous fixons l'orientation sur les courbes $ \beta_{i}. $
\begin{rem} \label{rem:bigone}
Si $  u\in \mathcal{M} \left( x,y \right)$ est un bigone de signe positif, alors $ Im(x)<Im(y) $ o\`u $ Im(x) $ d\'esigne la partie imaginaire de $ x $ vu comme un nombre complexe, voir la Figure \ref{fig:signe-bigone}. Si $ V_{x},V_{y} $ sont deux voisinages de $ x $ et $ y $ respectivement, alors $ V_{x}\cap u $ et $ V_{y}\cap u $ sont des coins convexes qui se trouvent respectivement, dans le quatri\`eme et le premier quadrant.\\
Si le signe de $ u $ est n\'egatif, alors $ Im(x)>Im(y) $ et $ V_{x}\cap u , V_{y}\cap u $  se trouvent respectivement, dans le deuxi\`eme et troisi\`eme quadrant.
\end{rem} 
\paragraph*{}
Soit $u \in \mathcal{M}(x, y)$ vu comme une application. Nous appellons nombre de multiplicit\'e par rapport \`a $u  $  d'un point $ w\in D_{n} $ le nombre, not\'e $ n_{u}(w), $ de points dans l'image r\'eciproque de $ w $ par $ u. $ 
\begin{displaymath}
 n_{u}(w)= \sharp u^{-1}(w).
\end{displaymath}
Le comportement du nombre de multiplicit\'e pr\`es des arcs $ \alpha_{xy} $ et $ \beta_{xy} $ est comme suit. Si nous traversons $ \alpha_{xy} $ ou $ \beta_{xy} $ le nombre de multiplicit\'e fait une variation de un.
 Pr\`es des points $ x $ et $ y, $ le nombre de multiplicit\'e est \'egal \`a $ 1 $ \`a l'int\'erieur du bigone et $ 0 $ ailleurs, comme le prouve le lemme suivant.
\begin{proof}
\begin{figure}[htbp] 
\psfrag{alphaxy}{$ \alpha_{xy} $}
\psfrag{betaxy}{$ \beta_{xy} $}
\psfrag{Vx}{$ V_{x} $}
\psfrag{Vy}{$ V_{y} $}
\psfrag{gamma}{$ \gamma $}
\psfrag{gamma+}{$ \gamma_{+} $}
\psfrag{-i}{$ -i $}
\psfrag{+i}{$ +i $}
\centerline{\includegraphics[width=5cm]{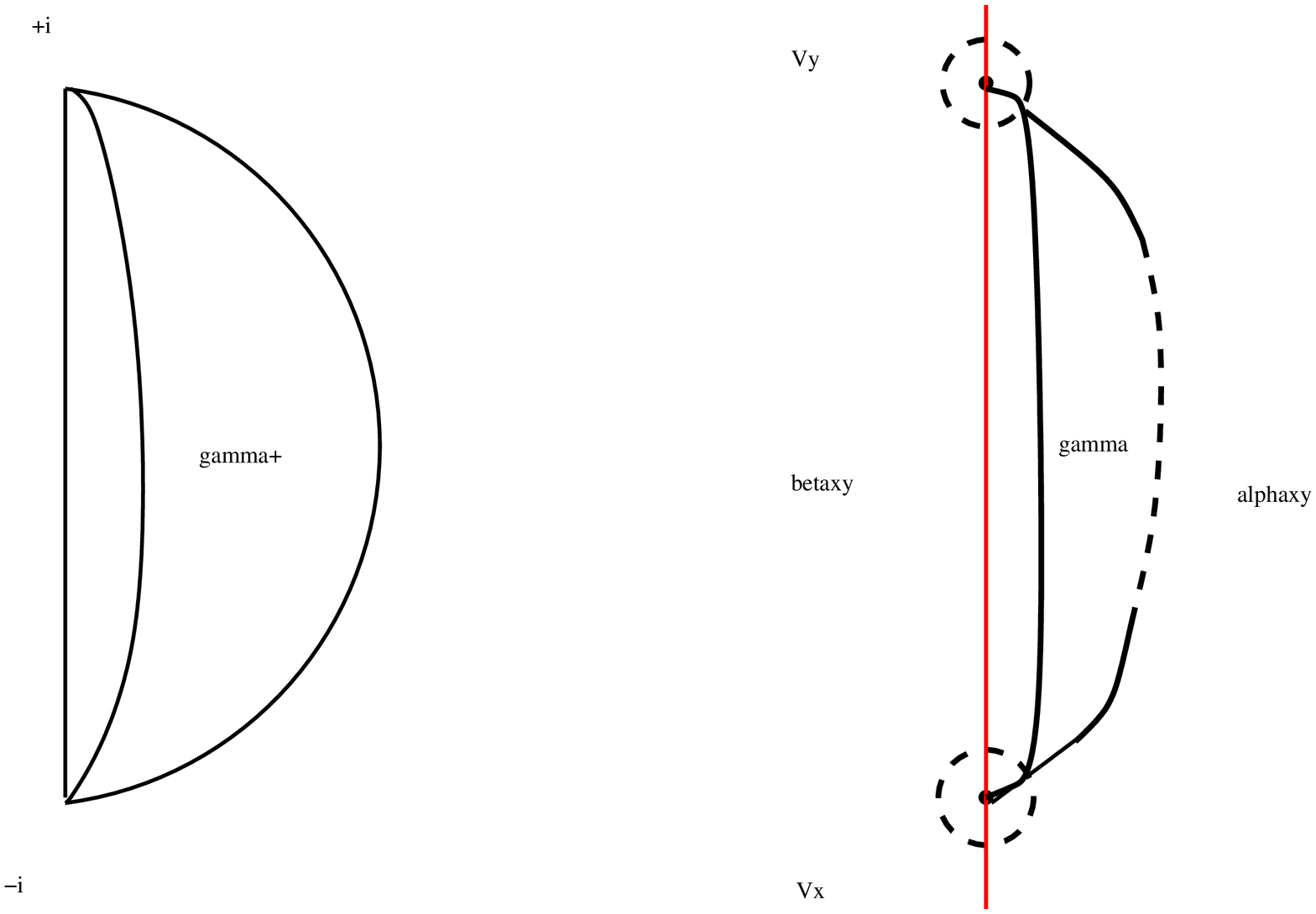}}
\caption{}\label{fig:multiplicite}
\end{figure}
\begin{lem} \label{lem:multiplicite-extremite}
Soit $ u $ un bigone de $ x $ vers $ y. $ Alors $ n_{u}(x)=n_{u}(y)=1. $
\end{lem}
Soit $ \gamma $ un arc plong\'e dans $ D_{n} $ d'extr\'emit\'es $ x $ et $ y $ telle que $ \gamma \cup \beta_{xy} $ borde un petit disque, not\'e $ \mathcal{D}', $ dans $ D_{n}, $ voir la Figure \ref{fig:multiplicite}. Nous voyons $ \mathcal{D}' $ comme l'image de la surface bord\'e par $ \left[ -i,i\right] \cup \gamma_{+} $ dans le demi disque $ \mathbb{D}$ par un plongement $ u', $ o\`u $ \gamma_{+} $ est un arc simple plong\'e dans $ \mathbb{D} $ envoy\'e sur $ \gamma $ par $ u' $ . Soit $ V_{x} $ et $ V_{y} $ deux voisinages de $ x $ et $ y $ respectivement, nous pouvons choisir ces voisinages suffisament petits de telle sorte que les morceaux de courbes de $ \gamma $ et $ \alpha_{xy} $ coincident \`a l'int\'erieur de ces voisinages i.e  
\begin{eqnarray}
\gamma \cap V_{x}= \alpha_{xy} \cap V_{x},
\gamma \cap V_{y}= \alpha_{xy} \cap V_{y}.
\end{eqnarray}
Alors les arcs obtenus de $ \gamma $ et $ \alpha_{xy} $ dans la surface $ D_{n} $ priv\'ee de $ \gamma\cap V_{x} $ et $ \gamma\cap V_{y} $ sont homotopes \`a extr\'emit\'es fixes (ici l'homotopie vient du fait que $ \mathbb{D}\cap S^{1} $ et $ \gamma_{+} $ sont homotopes dans $ \mathbb{D} $). Comme $ u' $ est un plongement, alors d'apr\`es (Th\'eor\`eme 3.1, \cite{EPSTEIN}), ces arcs sont isotopes \`a extr\'emit\'es fixes. Donc nous pouvons r\'eduire le nombre de multiplicit\'e \`a l'ext\'erieur des deux voisinages $ V_{x} $ et $ V_{y} $ sans le changer \`a l'int\'erieur. D'o\`u
\begin{eqnarray} 
n_{u}(w)=n_{u'}(w)=1,\text{ } \forall w\in V_{x}\cap \mathcal{D}'
\end{eqnarray}
\end{proof}

\section{Complexe de cochaines} \label{sec:complexedechaine}
Etant donn\'ee une tresse \`a $ n $ brins $ \sigma. $ Choisissons $ h_{\sigma} \in \text{Diff}^{+}(D,P_{n})$ un repr\'esentant de $ \sigma. $ Etant donn\'es  $ \alpha_{j},\beta_{i} \left( i,j\in \left[ \left[ 1,...,n-1\right] \right]  \right) $ deux arcs comme dans \mbox{la Figure \ref{fig:base},} nous associons \`a $ \left( \beta_{i},h_{\sigma}(\alpha_{j}) \right) $  un complexe de cocha\^ines $\left(  CF^{l,k}\left( \beta_{i},h_{\sigma}(\alpha_{j}) \right),\partial^{l,k}\right)  , l,k\in \Z. $  
 $ CF^{l,k}\left( \beta_{i},h_{\sigma}(\alpha_{j}) \right)  $ est le groupe ab\'elien libre engendr\'e par les points d'intersections de $ \beta_{i} $ avec $ h_{\sigma}(\alpha_{j}), $ o\`u  $ k,l $ sont respectivement les degr\'es d'Alexander et de Maslov que nous allons les d\'efinir.\\
 Nous munissons l'ensemble des g\'en\'erateurs d'un degr\'e d'Alexander \mbox{$A(x)=k, $} o\`u $ k $ est l'exposant paru dans la formule (\ref{eq:Ialgebrique}) section 2.
\paragraph*{}
Rappellons qu'une \textit{Fourchette} dans $ (D,P_{n}) $ est un arbre $ F $ plong\'e dans $ D $ form\'e de trois ar\^etes et quatre sommets $ d_{0},p_{i},p_{j} $ et $ z $ tels que $ F\cap \partial D=\left\lbrace d_{0}\right\rbrace , F\cap P_{n}=\left\lbrace p_{i},p_{j}\right\rbrace , $ et les trois ar\^etes ont $ z $ comme sommet commun. L'ar\^ete $ M(F) $ qui joint les deux sommets $ d_{0} $ et $ z $ s'appelle le \textit{manche} de $ F. $ L'union des deux autres ar\^etes est un arc plong\'e \`a extr\'emit\'es $ \left\lbrace p_{i},p_{j}\right\rbrace , $ appel\'e les \textit{dents} de $ F. $ Cette notion de fourchette a \'et\'e utilis\'e ind\'ependament par Dan Krammer en 2000, puis par Stephen Bigelow \cite{BIGELOW01} pour montrer que le groupe de tresses $ B_{n} $ est lin\'eaire.
\paragraph*{}
Pour tout $ j=1$,...,$n-1, $ notons par $ F_{j} $ la fourchette dans $ (D,P_{n}) $ qui a comme dents l'arc $ \alpha_{j} $ et $\partial M(F_{j})=\left\lbrace d_{0},z_{j}\right\rbrace  $  o\`u $ z_{j}$ est un point de $ \alpha_{j}. $ Nous identifions le point $ d_{0} $ avec l'un des bords des courbes $ \beta_{i}. $
\paragraph*{}
\begin{figure}[htbp] 
\psfrag{p1}{$ p_{1} $}
\psfrag{p3}{$ p_{3} $}
\psfrag{p4}{$ p_{4} $}
\psfrag{p2}{$ p_{2} $}
\psfrag{d0}{$ d_{0} $}
\psfrag{z2}{$ z_{2} $}
\psfrag{hz}{$ h_{\sigma_{2}}(z_{2}) $}
\psfrag{gama}{$ \gamma $}
\psfrag{x}{$ x $}
\centerline{\includegraphics[width=12cm,height=4cm]{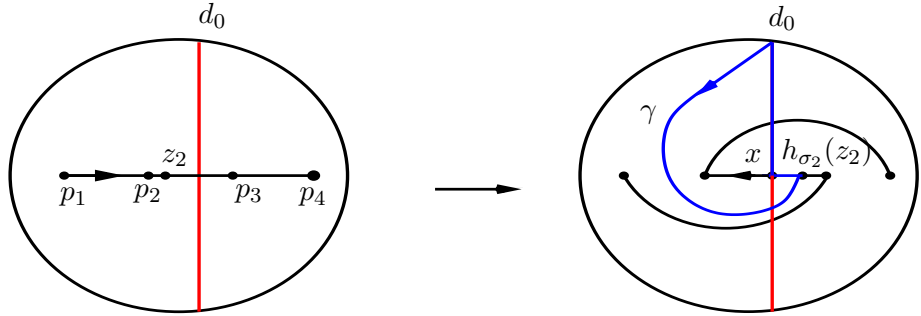}}
\caption{Exemple de chemin $ \gamma $}\label{fig:chemin}
\end{figure}

Etant donn\'ee un g\'en\'erateur $ x\in \beta_{i}\cap h_{\sigma}(\alpha_{j}). $ D\'efinissons un chemin plong\'e dans $ D_{n} $ param\'etr\'e par $ \gamma : \left[ 0,1\right] \rightarrow D_{n} $ comme suit:\\
Notons par $ \gamma_{1} $ l'arc qui va de $ d_{0} $ \`a $ h_{\sigma}(z_{j}) $ le long de $ h_{\sigma}(M(F_{j})). $\\
$ \gamma_{2} $ l'arc qui va de $ h_{\sigma}(z_{j}) $ vers $ x $ le long de $ h_{\sigma}(\alpha_{j}). $ \\
$ \gamma_{3} $ l'arc qui va de $ x $ \`a $ d_{0} $ le long de $ \beta_{i}. $\\
 Posons $ \gamma = \gamma_{1} \gamma_{2} \gamma_{3} ,$ voir la Figure \ref{fig:chemin}. Le chemin $ \gamma $ est $ C^{1} $ r\'egulier par morceaux. Nous d\'eformons l\'eg\'erement les courbes $ \beta_{i} $ et $ h_{\sigma}(\alpha_{j}) $ au voisinage des points $ x,h_{\sigma}(z_{j}) $ de sorte que les angles au points d'intersection seront droits ($ \pm\frac{\pi}{2} $).  Nous munissons l'ensemble des g\'en\'erateurs d'un degr\'e de Maslov $ \mu $ d\'efini par 
\begin{displaymath}
\mu\left( x \right) =\frac{1}{\pi}\int_{0}^{1}
d(arg(\gamma^{'}(t)))
\end{displaymath}
O\`u $ arg(\gamma^{'}(t)) $ d\'esigne l'argument du vecteur tangent $ \gamma^{'}(t). $ Autrement dit, $ \mu\left( x \right) $ est le nombre de demi-tours faits par le vecteur tangent $ \gamma^{'}(t). $ 
\paragraph*{}
La diff\'erentielle $ \partial : CF(\beta_{i},h_{\sigma}(\alpha_{j})) \rightarrow CF(\beta_{i},h_{\sigma}(\alpha_{j})) $ est d\'efinie sur un g\'en\'erateur $ x $ par
 \begin{equation}
\partial x = \sum_{\substack{y \in \beta_i \cap h_{\sigma}( \alpha_{j}) \\ u \in \mathcal{M}(x,y)}} \varepsilon \left( u \right) .y.
\end{equation}
o\`u $ \varepsilon\left( u\right)  $ est le signe du bigone $ u $ d\'efini dans la section 2. 
\begin{prop}
La diff\'erentielle $ \partial $ augmente  le degr\'e de Maslov par un et respecte le degr\'e d'Alexander.
\end{prop}
\begin{proof}
Soient $ x,y $ deux g\'en\'erateurs de $ CF^{l,k} \left( \beta_{i}, h_{\sigma}(\alpha_{j})\right).$ Si $ u$ est un bigone \'el\'ementaire de $ x $ vers $ y, $ alors le vecteur tangent fait un demi tour le long de $ \alpha_{xy}. $ Donc nous avons $ \mu (x,y) = \mu(y)-\mu (x) = 1. $ Si $ u $ n'est pas \'el\'ementaire. Alors nous pouvons nous ramener au premier cas par une isotopie qui \'elimine des bigones \'el\'ementaires. Mais 
  le nombre de demi-tours est invariant par isotopie, car un entier qui se d\'eforme d'une mani\`ere continue est forc\'ement constant. Donc $ \mu (x,y) =1. $ 
\paragraph*{}
Montrons que $ \partial $ respecte le degr\'e d'Alexander $ k. $ Supposons que $ A(x)=k $ et $ A(y)=k', $ alors 
 $ k-k' $ est le nombre de tours que l'on fait autour des perforations en parcourant la courbe qui va de $ x $ \`a $ y $ le long de $ h_{\sigma}(\alpha_{j}) $ et de $ y $ \`a $ x $ le long de $ \beta_{i}. $ Mais cette courbe borde un bigone. Donc $ k-k'=0. $ D'o\`u $ k=k'. $
\end{proof}
\begin{figure}[htbp] 
\psfrag{x}{$ x $}
\psfrag{y}{$ y $}
\psfrag{y0}{$ y_{0} $}
\psfrag{z}{$ z $}
\psfrag{u}{$ u $}
\psfrag{u0}{$ u_{0} $}
\psfrag{v}{$ v $}
\psfrag{v0}{$ v_{0} $}
\centerline{\includegraphics[width=8cm, height=3cm]{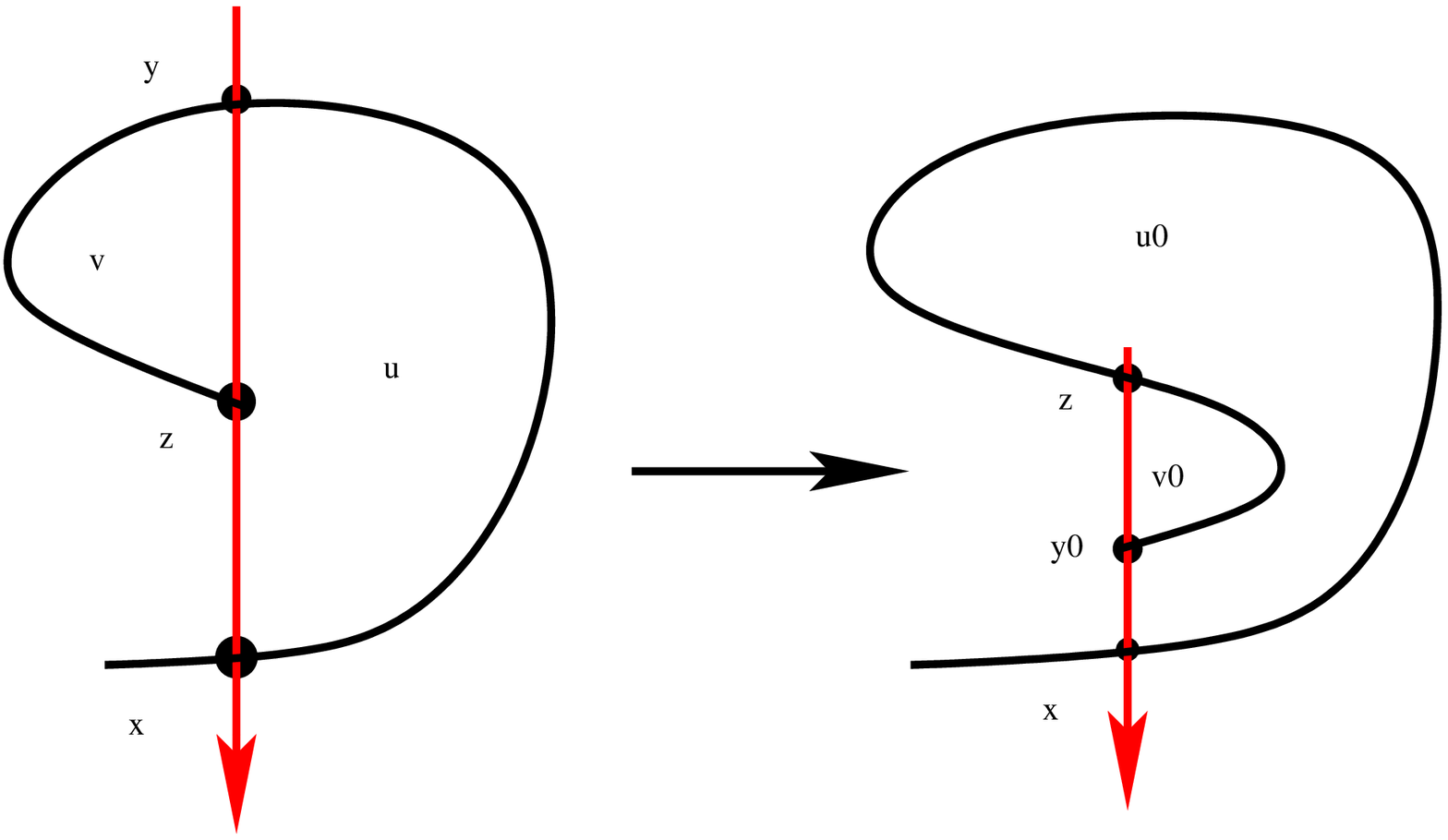}}
\caption{}\label{fig:bord}
\end{figure}
\begin{figure}[htbp] 
\psfrag{x}{$ x $}
\psfrag{y}{$ y $}
\psfrag{y0}{$ y_{0} $}
\psfrag{z}{$ z $}
\psfrag{u}{$ u $}
\psfrag{u0}{$ u_{0} $}
\psfrag{v}{$ v $}
\psfrag{v0}{$ v_{0} $}
\centerline{\includegraphics[width=8cm, height=3cm]{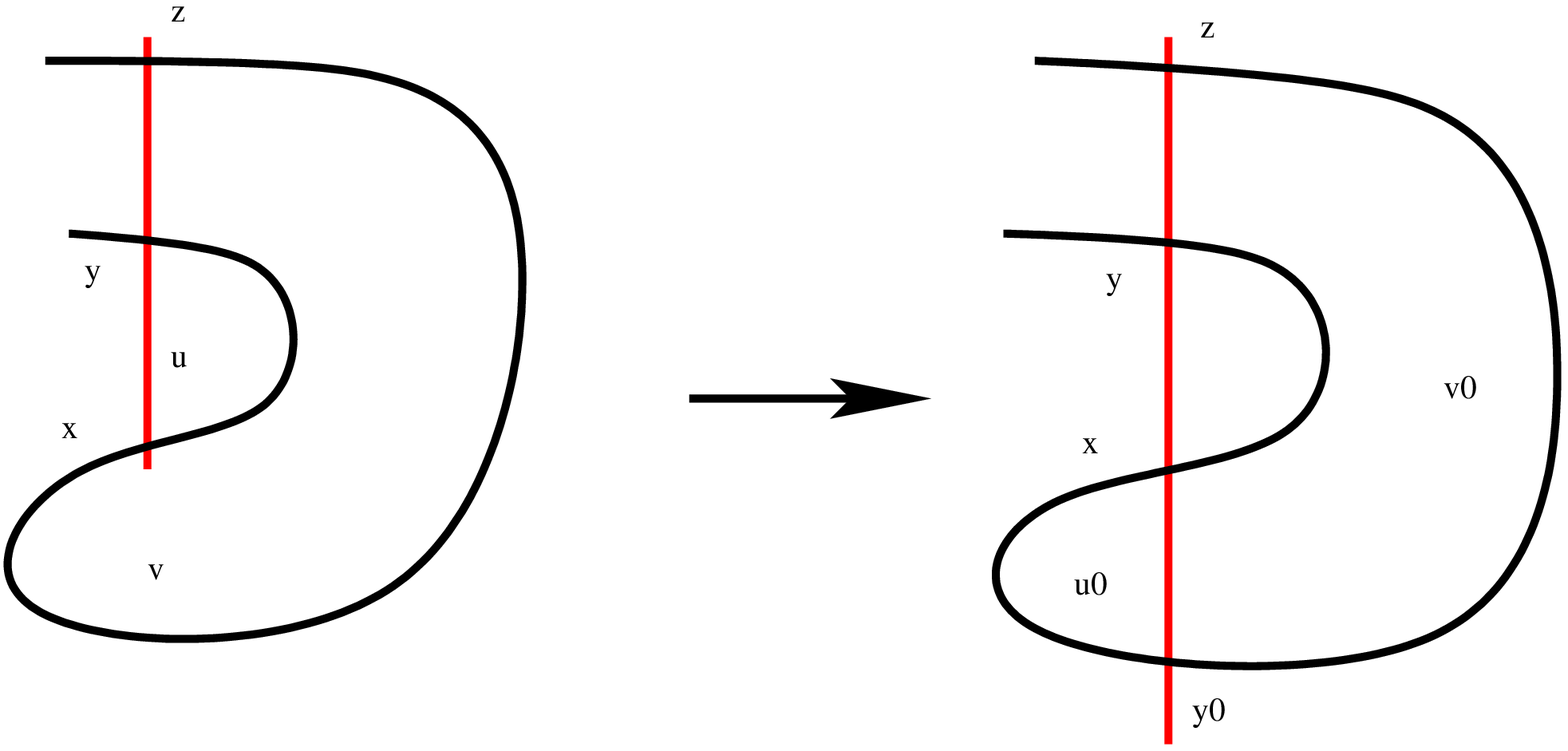}}
\caption{}\label{fig:bord1}
\end{figure}
\begin{prop} \label{pro:bord}
L'homomorphisme $ \partial $ est une diff\'erentielle. i.e $ \partial\circ \partial =0. $
\end{prop}
\begin{proof}

Soit $ x \in \beta_{i} \cap h_{\sigma}(\alpha_{j}),$ vu comme un g\'en\'erateur de $CF\left( \beta_{i},h_{\sigma}(\alpha_{j}) \right).$ Donc 
\begin{equation} \label{eq:bord}
\partial\circ \partial (x)= \sum_{\substack{y \in \beta_i \cap h_{\sigma}(\alpha_{j}) \\ u \in \mathcal{M}(x,y)}} \sum_{\substack{z \in \beta_i \cap h_{\sigma}(\alpha_{j}) \\ v \in \mathcal{M}(y,z)}}   \varepsilon \left( u \right). \varepsilon \left( v \right).z.
\end{equation}                    
Nous montrons que le coefficient de chaque g\'en\'erateur $ z $ s'annule. Supposons qu'il existe deux bigones $ u\in \mathcal{M}(x,y) $ et $ v\in \mathcal{M}(y,z). $ D'apr\`es la Remarque \ref{rem:bigone}, il y a plusieurs possibilit\'es pour que les bigones $ u $ et $ v $ occupent les quadrants pr\'es de $ y. $ Nous traitons deux cas, le m\^eme argument sera utilis\'e pour les autre cas.
\paragraph*{}
Dans le premier, supposons que $ u $ occupe le premier quadrant et $ v $ occupe le deuxi\`eme,  Les courbes $ \beta_{xy} $ et $ \beta_{yz} $ se chevauchent pr\`es de $ y. $ Donc l'un des courbes doit contenir l'autre, pour des raisons de muliplicit\'e locale. Supposons que $ \beta_{yz}\subset \beta_{xy}, $ l'autre cas se fait de la m\^eme mani\`ere. Comme l'int\'erieur de $ u\cup v $ ne contient pas des perforations et la courbe $ h_{\sigma}(\alpha_{j}) $ est simple, alors la courbe $ h_{\sigma}(\alpha_{j}) $ doit couper transversallement $ \beta_{xy}, $ le bord du bigone $ u $ qui se trouve sur la courbe $ \beta_{i} ,$ en un premier point $ y_{0} $ lorsque l'on parcourt $ h_{\sigma}(\alpha_{j}) $ en partant de $ x. $\\
Maintenant, l'intersection de la courbe $ h_{\sigma}(\alpha_{j}) $ avec $ \beta_{xy} $ en $ y_{0} $ va cr\'eer deux autres bigones $ u_{0}\in \mathcal{M}(x,y_{0}) $ et $ v_{0}\in \mathcal{M}(y_{0},z) ,$ voir le Figure \ref{fig:bord}. Les bigones $ u $ et $ u_{0} $ sont de m\^eme signe or les bigones $ v $ et $ v_{0} $ ont des signes oppos\'es i.e les coefficients $ \varepsilon(u).\varepsilon(v),\varepsilon(u_{0}).\varepsilon(v_{0}) $ de $ z $ ont des signes oppos\'es.
\paragraph*{}
Dans le deuxi\`eme cas, supposons que le bigone $ u $ occupe le premier quadrant et $ v $ occupe le quatri\`eme. Les courbes $ \beta_{xy} $ et $ \beta_{yz} $ s'intersectent qu'aux point $ y, $ voir la figure \ref{fig:bord1}. Dans ce cas la courbe $ \alpha_{yz} $ doit contenir la courbe $ \alpha_{xy}. $ alors la courbe $ h_{\sigma}(\alpha_{j}) $ doit couper transversallement la courbe $ \beta_{i} $ en un premier point $ y_{0} $ avec $ Im(y_{0})<Im(x) $ lorsque l'on parcourt $ h_{\sigma}(\alpha_{j}) $ en partant de $ y $ en direction de $ z. $ Maintenant, l'intersection de la courbe $ h_{\sigma}(\alpha_{j}) $ avec $ \beta_{i} $ en $ y_{0} $ va cr\'eer deux autres bigones $ u_{0}\in \mathcal{M}(x,y_{0}) $ et $ v_{0}\in \mathcal{M}(y_{0},z) ,$ voir la Figure \ref{fig:bord1}. Les bigones $ u $ et $ u_{0} $ ont des signes oppos\'es or les bigones $ v $ et $ v_{0} $ sont de m\^eme signe  i.e les coefficients $ \varepsilon(u).\varepsilon(v),\varepsilon(u_{0}).\varepsilon(v_{0}) $ de $ z $ ont des signes oppos\'es.
\paragraph*{}  
Ceci signifie que les coefficients dans la formule (\ref{eq:bord}) s'annulent par paires. D'o\`u $ \partial\circ \partial (x)= 0. $
\paragraph*{}
Montrons que le $ y_{0} $ est unique dans le sens: si $ y_{1} $ un autre point d'intersection tels qu'ils existent un bigone de $  x$ \`a $y_{1}$ et un bigone de $y_{1}$ \`a $z$ alors $ y_{0}=y_{1}. $ Nous traitons que le premier cas, le m\^eme raisonement se fait pour le deuxi\`eme. Supposons qu'il existe un point $ y_{1} $ telle que $  \mathcal{M}(x,y_{1})\neq \Phi$ alors l'arc $ \alpha_{y_{1}z} $ doit passer par le point $ y_{0} $ lorsque nous parcourons $ h_{\sigma}(\alpha_{j}) $ en partant de $ y_{1}. $ Dans ce cas le point $ y_{1} $ doit \^etre sur l'arc $ \beta_{y_{0}z}, $ mais il peut pas y avoir de bigone de $ x $ vers les points de $ \beta_{y_{0}z}-\{y_{0}\} $ pour la raison de multiplicit\'e locale. Donc $ y_{1}=y_{0}. $

\end{proof}
Nous d\'efinissons $ HF^{l}\left( \beta_{i},h_{\sigma}(\alpha_{j}) \right)  $ comme \'etant le $ l $-i\`eme
groupe de cohomologie du complexe de cocha\^ines $ CF\left( \beta_{i},h_{\sigma}(\alpha_{j}) \right) . $ C'est un $ \Z $-module gradu\'e de type fini. Nous d\'efinissons $ HF^{l,k}(\beta_{i},h_{\sigma}(\alpha_{j})) $ comme \'etant le $ l $-i\`eme groupe de cohomologie du sous complexe de degr\'e $ k $ de $ CF\left( \beta_{i},h_{\sigma}(\alpha_{j}) \right) . $ Donc nous avons la d\'ecomposition
\begin{eqnarray}
HF^{l}(\beta_{i},h_{\sigma}(\alpha_{j}))=\oplus_{k\in \Z}HF^{l,k}(\beta_{i},h_{\sigma}(\alpha_{j})).
\end{eqnarray}
Notons par $ HF^{l,k}(\sigma,i,j) $ la classe d'isomorphisme de $ HF^{l,k}(\beta_{i},h_{\sigma}(\alpha_{j})). $
\begin{thm} \label{th:invariance}
Si $ \sigma $ est une tresse \`a $ n $ brins. Pour chaque couple d'indice $ \left( i,j \right) , $ $ i,j \in \left[ \left[  1,...,n-1 \right]\right]  ,  $ la classe d'isomorphisme $ HF^{l,k}\left( \sigma,i,j \right)  $ ne d\'epend pas de la classe d'isotopie de la tresse. 
\end{thm}
La d\'emonstration de ce th\'eor\`eme est report\'ee \`a la section \ref{sec:invariance}. Nous montrons maintenant  que le coefficient de la matrice de Burau r\'eduite est \'egal \`a la caract\'eristique d'Euler gradu\'ee de cette cohomologie. Pour cela nous avons besoin du lemme suivant:
\begin{lem} \label{lem:signe}
Soit $ x $ un point dans $ \beta_{i} \cap h_{\sigma}( \alpha _{j}). $ Alors $ \mu (x) $ est pair si et seulement si $ sign (x) = +1, $ o\`u $ sign (x) $ est le signe d'intersection des courbes orient\'ees $ \beta_{i} $ et $ h_{\sigma}(\alpha_{j}) $ au point $ x. $
\end{lem}
\begin{proof}
Supposons que $ \beta_{i} $ et $ h_{\sigma}( \alpha_{j}) $ se coupent aux points $ x $ et $ h_{\sigma}(z_{j}) $ \`a angle droit $ \left( \pm \frac{\pi}{2}\right) . $ Soit $ \tilde{\gamma} $ une courbe simple donn\'ee comme dans la Figure \ref{fig:maslov}. 
Nous lissons la courbe $ \gamma $ aux points $ x $ et $ h_{\sigma}(z_{j}). $ Nous obtenons $ \tilde{\gamma}\gamma  $  une courbe simple ferm\'ee. Nous distinguons deux cas.\\
$ \mathbf{Cas 1} : $ Si le signe d'intersection de $ \beta_{i} $ avec $ h_{\sigma}( \alpha_{j}) $ au point $ x $ est positif, alors 
\begin{displaymath}
\int d(arg((\tilde{\gamma}\gamma)' (t))) =  \int d(arg(\gamma' (t))) + \int d(arg(\tilde{\gamma}' (t))) + \pi,
\end{displaymath}
\begin{displaymath}
2n\pi = \int d(arg(\gamma' (t))) + 2\pi,
\end{displaymath}
o\`u $ n $ est le degr\'e de l'application $ \tilde{\gamma}\gamma : S^{1} \rightarrow S^{1}. $ Alors $ \mu(x) = 2n-2. $\\
$ \mathbf{ Cas 2 :} $ Si le signe d'intersection de $ \beta_{i} $ avec $ h_{\sigma}( \alpha_{j}) $ au point $ x $ est n\'egatif, alors
\begin{displaymath}
\int d(arg((\tilde{\gamma}\gamma)' (t))) =  \int d(arg(\gamma' (t))) + \int d(arg(\tilde{\gamma}' (t))).
\end{displaymath}
Donc $ \mu(x) = 2n-1. $\\
\begin{figure}[htbp] 
\psfrag{c1}{\bf{Cas 1}}
\psfrag{c2}{\bf{Cas 2}}
\psfrag{p1}{$ p_{1} $}
\psfrag{p3}{$ p_{3} $}
\psfrag{p4}{$ p_{4} $}
\psfrag{p2}{$ p_{2} $}
\psfrag{d0}{$ d_{0} $}
\psfrag{z}{$ h_{\sigma}(z_{j}) $}
\psfrag{gammat}{$ \tilde{\gamma} $}
\psfrag{sigma1}{$ h_{\sigma}(\alpha_{1}) $}
\psfrag{beta2}{$ \beta_{2} $}
\centerline{\includegraphics[width=10cm,height=4cm]{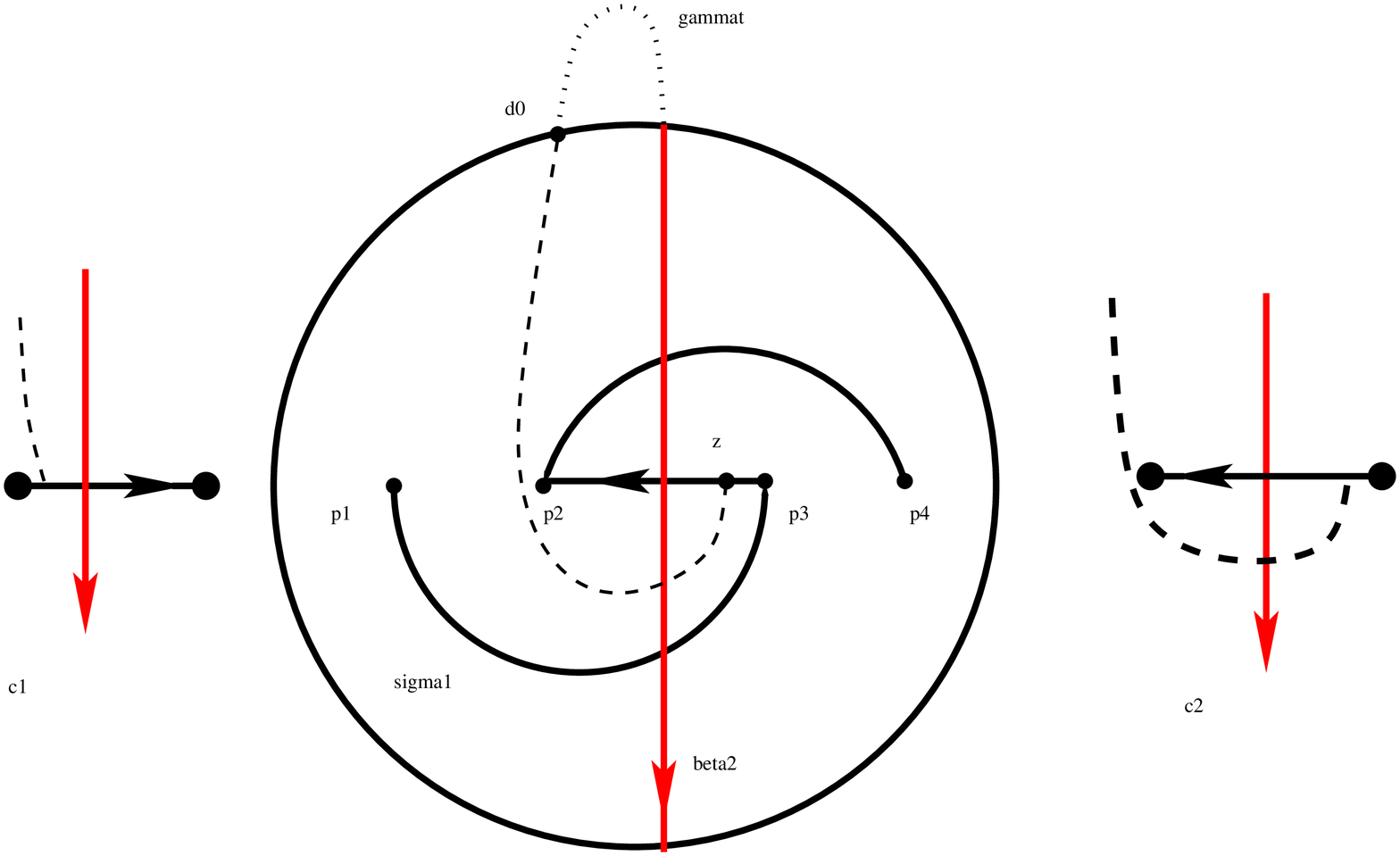}}
\caption{}\label{fig:maslov}
\end{figure}
\end{proof}
\begin{prop}
Soit $ \sigma $ une tresse \`a $ n $ brins. Pour chaque couple d'indice $ \left( i,j \right) , $ $ i,j \in \left[\left[  1,...,n-1 \right]\right] ,  $ il existe un complexe de cocha\^ines bigradu\'e $ \left( CF^{l,k}\left( \beta_{i},h_{\sigma}(\alpha_{j}) \right), \partial^{l,k} \right)_{l,k\in \Z}   $ tel que
\begin{equation} \label{eq:categorification}
b_{ij}= \sum _{l\in \Z}\left( -1\right) ^{l} . t^{k}.dim_{\Q}\left( HF^{l,k}(\beta_{i},h_{\sigma}(\alpha_{j}))\otimes_{\Z}\Q \right) .
\end{equation}
O\`u $ b_{ij} $ est le coefficient de la matrice de Burau $ \tilde{\rho}(\sigma) $ en position $ (i,j) .$ 
\end{prop}
\begin{proof}
C'est une cons\'equence directe du Lemme \ref{lem:signe} et la formule (\ref{eq:Ialgebrique})
 donn\'ee dans la section \ref{sec:courbe}.
\end{proof} 
\section{Invariance} \label{sec:invariance}
Dans cette section nous montrons que la cohomologie construite \`a la section \ref{sec:complexedechaine} ne d\'epend pas de la classe d'isotopie de la tresse. Nous utilisons pour cela l'\'elimination gaussienne \cite{BAR-NATAN07}, voir aussi \cite{BELIAKOVA10} pour une version plus \'el\'ementaire adapt\'ee \`a notre situation. Rappelons que le groupe de tresses $ B_{n} $ est le groupe des diff\'eomorphismes $ \text{Diff}^{+}(D,P_{n}) $ quotient\'e par la relation d'isotopie. Chaque diff\'eomorphisme (repr\'esente une tresse g\'eom\'etrique) agit sur le disque $ D_{n}. $ Si deux diff\'eomorphismes isotopes agiraient sur le disque \'epoint\'e, cette isotopie fait appara\^itre ou dispara\^itre des bigones (\'el\'ementaires).
\paragraph*{}
Soit $ \sigma $ une tresse \`a $ n $ brins. Dans cette section, $ h_{\sigma} $ et $ h'_{\sigma} $ sont deux diff\'eomorphismes isotopes par une isotopie qui fait appara\^itre un bigone \'el\'ementaire d'un point $ x $ vers un point $ y. $ Supposons que le diff\'eomorphisme $ h'_{\sigma} $ est celui qui fait appara\^itre un bigone \'el\'ementaire de plus. Notons par $ \mathcal{C} $ et $ \mathcal{C}' $ les deux complexes de cocha\^ines associ\'es \`a $ h_{\sigma} $ et $ h'_{\sigma} $ respectivement. Auparavant nous allons donner le lemme suivant:
\begin{lem} \label{lem:elimination}
(\textbf{Elimination Gaussienne \cite{BAR-NATAN07}}) Soit $ \mathcal{C}=(C^{*},\partial) $ un complexe de cocha\^ines sur $ \Z $ librement engendr\'e. Soit $ x\in C^{l} $ (resp. $ y\in C^{l+1} $) telle que $ C^{l}=\Z.x \oplus A $ (resp. $ C^{l+1}=\Z.y\oplus B $). Si $\varphi:\Z.x \rightarrow \Z.y $ est un isomorphisme de $ \Z $-modules,  alors le segment de complexe de $ \mathcal{C} $
\begin{equation} \label{eq:BeforeGE}
\begin{CD}
  \cdots \left[ C^{l-1}\right] @>{\begin{pmatrix}\kappa \\ \zeta \end{pmatrix}}>> \begin{bmatrix} x \\ A \end{bmatrix} @>{\begin{pmatrix} \varphi & \delta \\ \lambda & \xi \end{pmatrix}}>> \begin{bmatrix} y \\ B \end{bmatrix} @>{\begin{pmatrix} \theta & \nu \end{pmatrix}}>> \left[ C^{l+2}\right] \cdots
\end{CD}
\end{equation}
est isomorphe au segment de complexe de cocha\^ines suivant
\begin{equation} \label{eq:duringGE}
\begin{CD}
  \cdots \left[ C^{l-1}\right] @>{\begin{pmatrix}0 \\ \zeta \end{pmatrix}}>> \begin{bmatrix} x \\ A \end{bmatrix} @>{\begin{pmatrix} \varphi & 0 \\ 0 & \xi-\lambda\varphi^{-1}\delta \end{pmatrix}}>> \begin{bmatrix} y \\ B \end{bmatrix} @>{\begin{pmatrix} 0 & \nu \end{pmatrix}}>> \left[ C^{l+2}\right] \cdots
\end{CD}
\end{equation}
ces deux complexes sont homotopiquement \'equivalent au segment de complexe 
\begin{equation} \label{eq:afterGE}
\begin{CD}
  \cdots \left[ C^{l-1} \right] @>{\left(  \zeta \right)  }>> \left[  A \right]  @>{\left( \xi-\lambda\varphi^{-1}\delta \right) }>> \left[ B \right]   @>{\left(  \nu \right)  }>> \left[ C^{l+2} \right] \cdots
\end{CD}
\end{equation} 
Ici on utilise la notation matricielle pour la diff\'erentielle $ \partial. $
\end{lem}
\begin{figure}[htbp] 
\psfrag{x}{$ x $}
\psfrag{y}{$ y $}
\psfrag{t}{$ t $}
\psfrag{z}{$ z $}
\centerline{\includegraphics[width=8cm ,height=3cm]{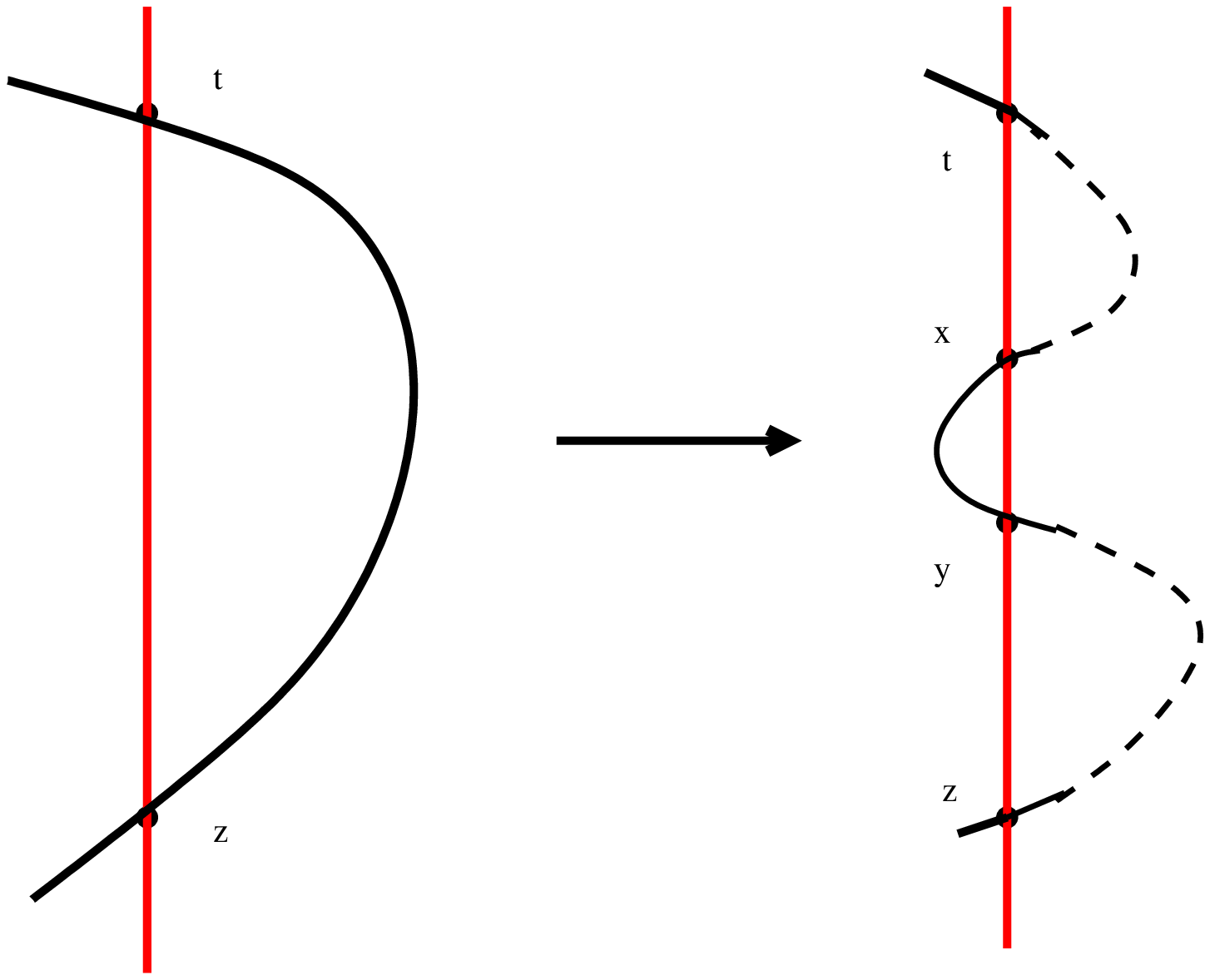}}
\caption{}\label{fig:3invariance}
\end{figure}
\begin{proof} [\bf{Preuve du th\'eor\`eme \ref{th:invariance}}].
Soient $ h_{\sigma} $ et $ h'_{\sigma} $ deux diff\'eomorphismes (repr\'esentant deux tresses g\'eom\'etriques) isotopes. Nous passons de $ h_{\sigma} $ \`a $ h'_{\sigma} $ par une isotopie qui fait appara\^itre un nombre $ n $ de bigones \'el\'ementaires. Par raissonons par induction sur le nombre de bigones.\\
\begin{figure}[htbp] 
\centerline{\includegraphics[width=10cm ,height=4cm]{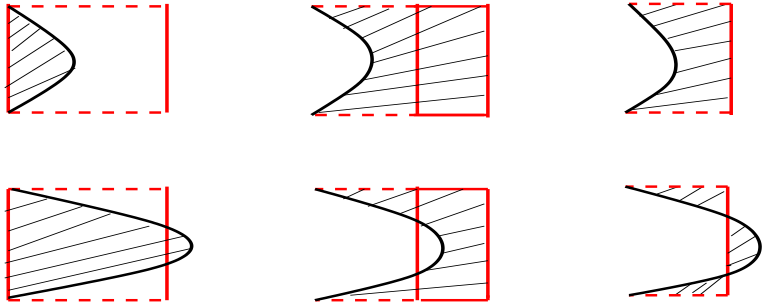}}
\caption{Bijection entre les bigones qui bordent l'arc $h_{\sigma}( \alpha_{j}) $ \`a l'endroit o\`u on \'elimine le bigone \'el\'ementaire, les hachures correspondent \`a une multiplicit\'e sup\'erieur de un. } \label{fig:bijection-bigone}
\end{figure}
Si $ n=1, $ on passe de $ h_{\sigma} $ \`a $ h'_{\sigma} $ par une isotopie qui ne fait appara\^itre qu'un bigone \'el\'ementaire d'un point $ x $ vers un point $ y $. Consid\'erons le complexe de cocha\^ines $ \mathcal{C}=(CF^{*,*}(\beta_{i},h_{\sigma}(\alpha_{j})),\partial) $ associ\'e au diff\'eomorphisme $ h_{\sigma} $
\begin{equation} \label{eq:avantisotopie}
\begin{CD}
  \cdots \left[ C^{l-1} \right] @>{\left( \partial \right)  }>> \left[  C^{l} \right]  @>{\left( \partial \right) }>> \left[ C^{l+1} \right]   @>{\left(  \partial \right)  }>> \left[ C^{l+2} \right] \cdots
\end{CD}
\end{equation}
Le complexe de cocha\^ines $ \mathcal{C}'=(CF(\beta_{i},h'_{\sigma}(\alpha_{j})),\partial') $ associ\'e au diff\'eomorphisme $ h'_{\sigma} $ est de la forme
\begin{equation} \label{eq:apresisotopie}
\begin{CD}
  \cdots \left[ C^{l-1}\right] @>{\begin{pmatrix}\kappa \\ \zeta \end{pmatrix}}>> \begin{bmatrix} x \\ C^{l} \end{bmatrix} @>{\begin{pmatrix} \varphi & \delta \\ \lambda & \xi \end{pmatrix}}>> \begin{bmatrix} y \\ C^{l+1} \end{bmatrix} @>{\begin{pmatrix} \theta & \nu \end{pmatrix}}>> \left[ C^{l+2}\right] \cdots
\end{CD}
\end{equation}
o\`u les lettres grecques sont des matrices d'homomorphismes de $ \Z $-modules. L'ensemble des g\'en\'erateurs de $ \mathcal{C}' $ est obtenu de celle de $ \mathcal{C} $ en rajoutant les points $ x $ et $ y. $ L'isomorphisme $ \varphi $ est d\'efini par  $ \varphi(x)=\varepsilon(u)y=\pm y  $ o\`u $ u\in \mathcal{M}(x,y), $ et $ \varphi(z)=0 $ pour $ z\neq x. $ Remarquons que $ \varphi^{-1}\circ \partial'(x)=x. $ Soit $ C'' $ le complexe de cocha\^ines obtenu apr\`es avoir appliqu\'e l'\'elimination gaussienne sur le complexe de cocha\^ines $ \mathcal{C}' $ associ\'e au diff\'eomorphisme $ h'_{\sigma}. $ La diff\'erentielle $ \partial''=\partial'-\partial'\circ \varphi^{-1}\circ \partial'. $ Montrons que $ \partial=\partial''. $ Soient $ z $ un g\'en\'erateur et $ t\in\partial z, $ donc il existe un bigone $ u\in \mathcal{M}(z,t). $ Si le bigone $ u $  a diparu lors du passage de $ h_{\sigma} $ \`a $ h'_{\sigma} $ par isotopie, alors la disparition de $ u $ fait appara\^itre un bigone $ u_{1} $ de $ z $ vers $ y $ et un bigone $ u_{2} $ de $ x $ vers $ t, $ voir la Figure \ref{fig:3invariance}. Si le bigone $ u $ n'a pas disparu, alors il y aura pas de bigones de $ z $ vers $ y $ ou de $ x $ vers $ t $, donc $ \varphi^{-1}\circ \partial'=0 $ ou $ \partial' \circ \varphi^{-1}=0. $ Donc les bigones compt\'es dans les diff\'erentielles $ \partial $ et $ \partial'' $ sont en bijection et ils viennent avec le m\^eme signe, voir la figure \ref{fig:bijection-bigone}.  
D'apr\`es le lemme d'\'elimination gaussienne les deux complexes de cocha\^ines $ \mathcal{C} $ et $ \mathcal{C}' $ sont homotopiquement \'equivalent.\\
Si $ n>1, $ nous r\'eappliquons l'\'elimination gaussienne. Nous obtenons que $ \mathcal{C} $ et $ \mathcal{C}' $ sont homotopiquement \'equivalent. Donc les groupes de cohomologie $ HF^{*,*}(\beta_{i},h_{\sigma}(\alpha_{j})) $ et $ HF^{*,*}(\beta_{i},h'_{\sigma}(\alpha_{j})) $ sont isomorphes.
\end{proof}
\section{Fid\'elit\'e}
Dans ce paragraphe, nous montrons que la cohomologie construite \`a la section \ref{sec:complexedechaine} d\'etecte la tresse triviale.
Nous fixons une collection de courbes $ \alpha_{0},...,\alpha_{n-1} $ comme dans la Figure \ref{fig:base}.
\begin{lem} \label{th:fidele}
\cite {KH-S} Si $ h $ un \'el\'ement de $ \text{Diff}^{+} \left( D, P_{n} \right)  $ satisfait $ h(\alpha_{i}) \simeq \alpha_{i} $ pour tout $ i, $ alors $ h\simeq id. $ 
\end{lem}
 Une courbe  $ c $  dans $ \left( D,P_{n} \right)  $ est dite \textbf{admissible} \cite{KH-S} si elle est l'image de $ \alpha_{i} $ par un \'el\'ement $ h $ de  $ \text{Diff}^{+} \left( D, P_{n} \right)  $ pour un certain $ i\in \left[ \left[  0,...,n-1 \right] \right]. $
\begin{figure}[htbp] 
\psfrag{alpha1}{ $\alpha_{1}$}
\psfrag{alphan-1}{$ \alpha_{n-1} $}
\psfrag{lg1}{$ l_{1}^{g} $}
\psfrag{lgn-2}{$ l_{n-2}^{g} $}
\psfrag{ld1}{$ l_{1}^{d} $}
\psfrag{ldn-2}{$ l_{n-2}^{d} $}
\centerline{\includegraphics[width=10cm ,height=4cm]{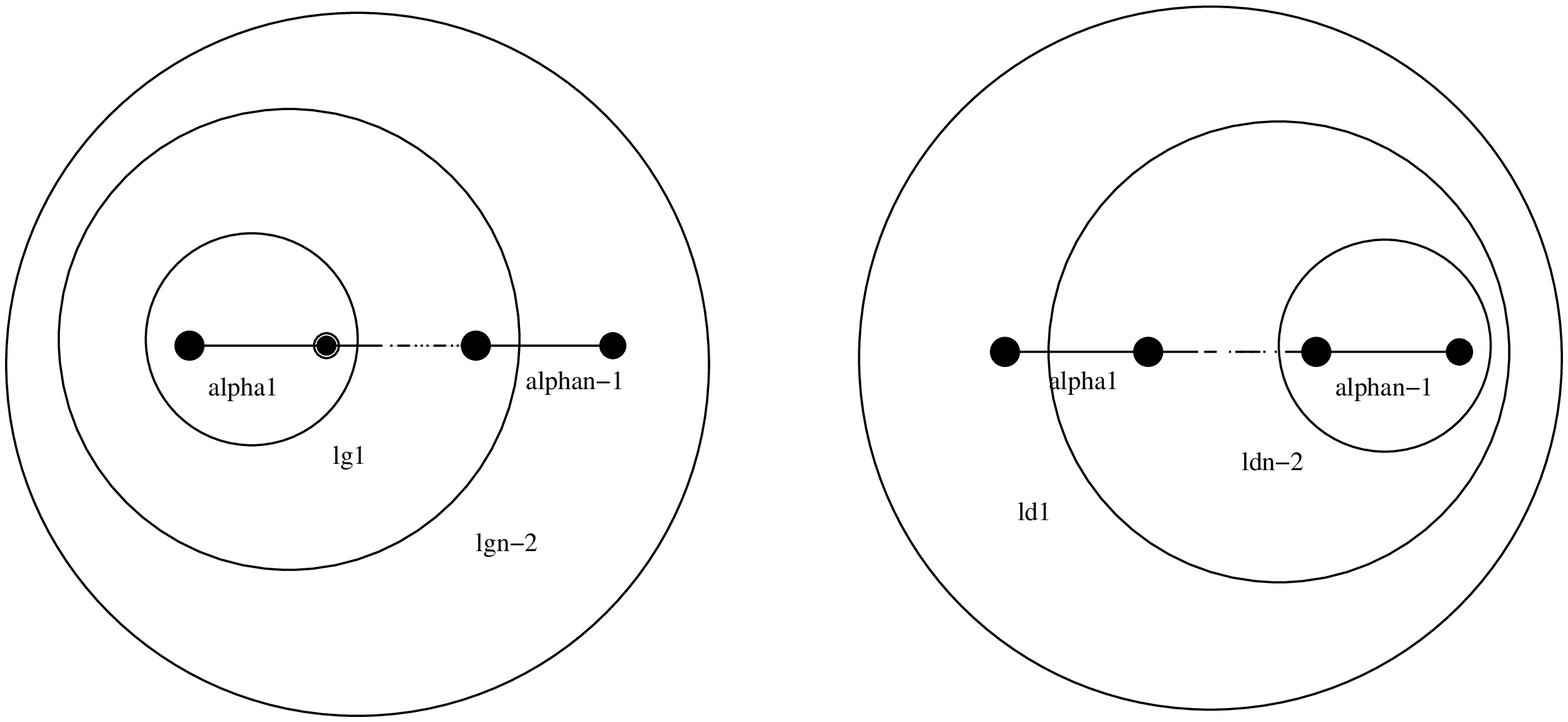}}
\caption{} \label{fig:isofidele}
\end{figure}
\begin{lem} \label{th:iso-fidele}
Soit $ c $ une courbe admissible dans $ (D,P_{n}). $ Supposons qu'il existe $ k\in \left[ \left[ 1,...,n-1 \right] \right] $ tel que $ i(\alpha_{j},c) = i(\alpha_{j},\alpha_{k}) $ pour tout $ j= 1$,...,$n-1. $ Alors on a
\[
c \simeq
\begin{cases}
 \alpha_1 \text{ ou } \tau_{d,1}^{\pm1}(\alpha_1) & \text{ si } k = 1,\\
 \alpha_k \text{ ou } \tau_{d,k}^{\pm1}(\alpha_k), \tau_{g,k}^{\pm1}(\alpha_{k}),\tau_{d,k} \tau_{g,k}^{-1}(\alpha_{k}) \text{ ou } \tau_{d,k}^{-1} \tau_{g,k}(\alpha_{k})  & \text{ si } 2 \leq k < n-2, \\
 \alpha_{n-1} \text{ ou } \tau_{g,n-1}^{\pm1}(\alpha_{n-1}) & \text{ si } k = n-1.
\end{cases}
\]
O\`u $ \tau _{d,1}$,...,$\tau _{d,n-2} $ sont les twists de Dehn positif le long des courbes ferm\'ees $ l_{1}^{d}$,...,$l_{n-2}^{d} $ et $ \tau _{g,1}$,...,$\tau _{g,n-2} $ sont les twists le long des courbes $ l_{1}^{g}$,...,$l_{n-2}^{g}, $ voir la Figure \ref{fig:isofidele}.
\end{lem}
\begin{proof}
\begin{figure}[htbp] 
\psfrag{p1}{ $p_{1}$}
\psfrag{pn}{$ p_{n} $}
\psfrag{neg}{$ c\simeq \tau^{-}_{d,1}(\alpha_{1}) $}
\psfrag{pos}{$ c\simeq \tau^{+}_{d,1}(\alpha_{1}) $}
\psfrag{cisob1}{$ c\simeq \alpha_{1} $}
\centerline{\includegraphics[width=12cm, height=6cm]{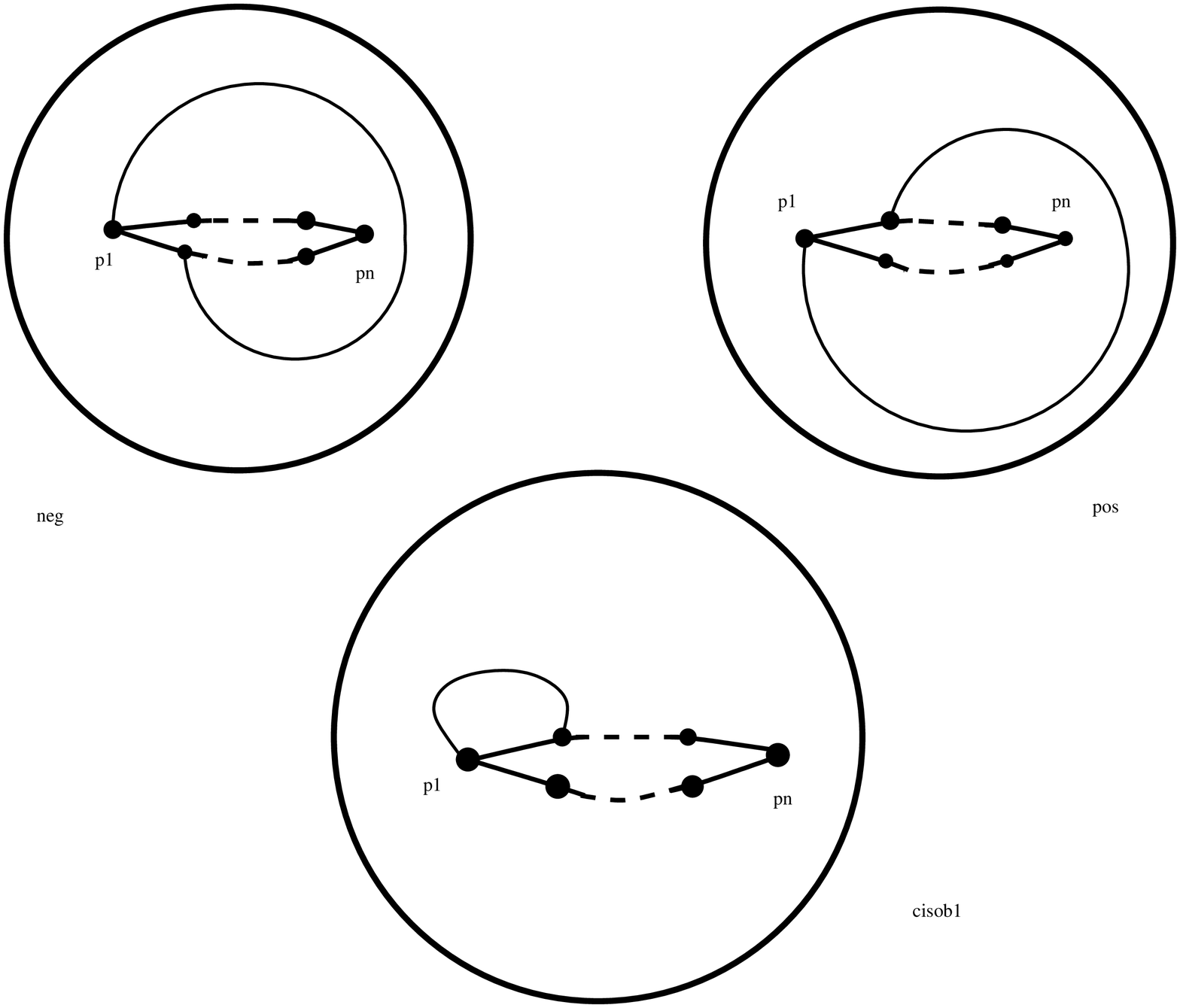}}
\caption{} \label{fig:isofideled}
\end{figure}
Montrons le r\'esultat pour le cas $ k=1, $ les autres cas se d\'emontrent d'une mani\`ere similaire. D'apr\`es l'hypoth\`ese, nous avons $ i(\alpha_{2},c)=i(\alpha_{2},\alpha_{1})=1. $ Donc les courbes $ c $ et $ \alpha_{1} $ ont un seul point commun qui est soit $ \alpha_{1}\cap \alpha_{2} $ ou soit $ \alpha_{2}\cap\alpha_{3}, $ mais le deuxi\`eme cas est impossible puisque $ i(\alpha_{3},c)=i(\alpha_{3},\alpha_{1})=0. $ D'autre part, nous avons $  i(\alpha_{1},c)=i(\alpha_{1},\alpha_{1})=2.$ Ceci signfie que les courbes $ c $ et $ \alpha_{1} $ ont les m\^emes extr\'emit\'es. Nous pouvons supposer que l'intersection de $ c $ avec tous les courbes $ \alpha_{j} $ est minimal. La courbe $ c $ ne rencontre  $ \alpha_{1}\cup ...\cup \alpha_{n-1} $ qu'aux points $ p_{1},p_{2}, $ encore par hypoth\`ese d'intersection g\'eom\'etrique. Maitenant nous coupons la surface $ D^{2} $ le long de $ \alpha_{1}\cup ...\cup \alpha_{n-1}. $ 
Nous consid\'erons la courbe $ c $ sur cette nouvelle surface, nous obtenons $ 3 $ cas possibles, 
voir la Figure \ref{fig:isofideled}. D'o\`u le r\'esultat.  
\end{proof}
\begin{thm}
Soit $ \sigma $ une tresse dans $ B_{n}. $ Supposons que pour tout couple d'indice $ \left( i,j \right), $ $ i,j \in \left[ \left[  1,...,n-1 \right] \right] , $ il existe un isomorphisme bigradu\'e de $ HF^{*,*} (\beta_{i},h_{\sigma}(\alpha_{j}) )$ vers $HF^{*,*} ( \beta_{i}, \alpha_{j} ). $ Alors $ \sigma$ est isotope \`a la tresse triviale.
\end{thm}
\begin{proof}
Supposons que $  HF^{*,*} (\beta_{i},h_{\sigma}(\alpha_{j}) ), HF^{*,*} ( \beta_{i},\alpha_{j} ) $ sont isomorphes. Alors pour tout $ i=1$,..,$n-1, $ $ HF^{*,*} (\beta_{i},h_{\sigma}(\alpha_{i} )$  est isomorphe au groupe ab\'elien libre engendr\'e par l'\'el\'ement $ x_{ii}, $ o\`u $ x_{ii} $ est le point d'intersection de $ \beta_{i} $ avec $ \alpha_{i}, $ et $ HF^{*,*} (\beta_{i},h_{\sigma}(\alpha_{j}) ) $ est isomorphe au groupe trivial si $ i\neq j. $ Alors
\begin{equation} \label{eq:triv-intersection}
i(\beta_{i},h_{\sigma}(\alpha_{j}))=
\begin{cases}
1 & \text{ si } i=j,\\
0 & \text{ si } i\neq j.
\end{cases}
\end{equation}
Donc $ \alpha_{j},h_{\sigma}(\alpha_{j}) $ ont les m\^emes ext\'emit\'es pour tout $ j=1$,..,$n-1. $ Ceci implique que 
\begin{displaymath}
i(\alpha_{k},h_{\sigma}(\alpha_{j}))=i(\alpha_{k},\alpha_{j}), \text{ pour tout } k=1,..,n-1.
\end{displaymath}
Appliquons le Lemme \ref{th:iso-fidele} pour $ c= h_{\sigma}(\alpha_{j}), $ Nous aurons
\begin{displaymath}
h_{\sigma}(\alpha_{j}) \simeq \alpha_{j} \text{ ou }
\tau_{d,j}^{\pm 1}(\alpha_j)  \text{ ou } \tau_{g,j}^{\pm1}(\alpha_{j}).
\end{displaymath}
Puisque $ i(\beta_{i},h_{\sigma}(\alpha_{j}))=0 $ pour tout $ i\neq j, $ alors $ h_{\sigma}(\alpha_{j}) \simeq \alpha_{j}, \forall j=1,..,n-1. $ D'apr\`es le Lemme \ref{th:fidele} on a $ h_{\sigma} \simeq (\tau_{0})^{m}, $ o\`u $ m\in \Z $ et $ \tau_{0} $ est le twist de Dehn positif le long du bord $ \partial D_{n}. $ Comme l'isomorphisme entre les groupes de cohomologie $  HF^{*,*} (\beta_{i},h_{\sigma}(\alpha_{j}) )$ et $HF^{*,*} ( \beta_{i}, \alpha_{j} )$ est bigradu\'e pour tout $ i,j\in \left[ \left[ 1,...,n-1\right] \right]  $, alors $ m=0. $ D'o\`u $ \sigma $ est la tresse triviale.
\end{proof}
\paragraph*{Conclusion et perspective}
Dans le travail de M. Khovanov et P. Seidel \cite{KH-S} la repr\'esentation par foncteur contient la notion de multiplicacit\'e qui n'apparait pas dans ce travail. M. Abouzaid \cite{ABOUZAID} a donn\'e une version compl\`ete de la cat\'egorie de Fukaya des surfaces ferm\'ees. Pouvons nous adapter le travail de Abouzaid pour decrire la structure multiplicative dans ce cadre ? 
\paragraph*{Remerciements}
Ce travail a \'et\'e r\'ealis\'e lorsque j'ai visit\'e le LMAM, Universit\'e de Bretagne Sud (France). Je tiens \`a remercier vivement mon directeur de th\`ese Christian Blanchet pour son soutien, pour ses encouragements et tous les conseils qu'il m'a prodigu\'e. Je remercie aussi Bertrant patureau-Mirand avec qui j'ai longuement discut\'e.

\end{document}